\documentclass[11pt,fullpage,leqno]{article}
\usepackage[english]{babel}
\usepackage{amsmath}
\numberwithin{equation}{section}
\usepackage{amsfonts}
\usepackage{amsmath}
\usepackage{epsfig}
\usepackage{graphicx}
\usepackage{color}
\usepackage[title,titletoc,toc]{appendix}

\def\R{\hbox{\bf R}}
\def\Z{\hbox{\bf Z}}
\def\1{\hbox{\bf 1}}

\makeatletter
\let\@fnsymbol\@arabic
\makeatother

\newcommand{\ba}{\begin{eqnarray}}
\newcommand{\ea}{\end{eqnarray}}

\newtheorem{theo}{\bf Theorem}[section]
\newtheorem{lem}[theo]{\bf Lemma}
\newtheorem{pro}[theo]{\bf Proposition}
\newtheorem{cor}[theo]{\bf Corollary}
\newtheorem{defi}[theo]{\bf Definition}
\newtheorem{rem}[theo]{\bf Remark}

\renewcommand{\R}{{\mathbb R}}

\renewcommand{\Z}{{\mathbb Z}}

\newenvironment{Proofc}[1]{\smallskip\par\noindent\textsc{#1}\quad}%
  {\hfill$\Box$\bigskip\par}

\nonstopmode

\begin{document}

\title{\bf Velocity diagram of traveling waves\\ for  discrete reaction-diffusion equations}

\author{M. Al Haj{}\footnote{Lebanese University, Faculty of Science (section 5), Nabatiye 1700 Lebanon.\newline E-mail: mohammad.alhaj.1@ul.edu.lb}
\and R. Monneau{}\footnote{CERMICS, Universit\'e Paris-Est, Ecole des Ponts ParisTech, 6-8 avenue Blaise Pascal, 77455 Marne-la-Vall\'ee Cedex 2, France; et 
CEREMADE, Universit\'e Paris-Dauphine, Place du Mar\'echal De Lattre De Tassigny, 75775 Paris Cedex 16, France.}}

\maketitle\thispagestyle{empty}


\noindent{\bf Abstract:} 
We consider a discrete version of reaction-diffusion equations. A typical example is the fully overdamped Frenkel-Kontorova model, where the velocity is proportional to the force. We also introduce an additional exterior force denoted by $\sigma$. For general discrete and fully nonlinear dynamics, we study traveling waves of velocity $c=c(\sigma)$ depending on the parameter $\sigma$. Under certain assumptions, we show properties of the velocity diagram $c(\sigma)$ for $\sigma\in [\sigma^-,\sigma^+]$. We show that the velocity $c$ is nondecreasing in $\sigma \in (\sigma^-,\sigma^+)$ in the bistable regime, with vertical branches $c\ge c^+$ for $\sigma=\sigma^+$ and $c\le c^-$ for $\sigma=\sigma^-$ in the monostable regime.\\



\section{Introduction}

\subsection{General motivation}

Our initial motivation was to study the classical fully overdamped Frenkel-Kontorova model,
which is a system of ordinary differential equations 
\begin{equation}\label{e1}
\frac{d X_{i}}{dt}=X_{i+1}-2X_{i}+X_{i-1}+f(X_{i})+\sigma \quad \mbox{with}\quad f(x)=-\beta\cos(2\pi x),\quad \beta>0
\end{equation}
where $X_{i}(t)\in\R$ denotes the position of a particle $i\in\Z$ at time
$t,$ $\displaystyle\frac{d X_{i}}{dt}$ is the velocity of this particle. Here $f$ is the force created by a $1$-periodic potential and $\sigma$ represents the constant driving force.
This kind of system can be, for instance, 
used as a model of  the motion of a dislocation defect in a crystal (see the book of Braun and Kivshar \cite{BK1}).
This motion is described by particular solutions of the form
\begin{equation}\label{e3}
X_{i}(t)=\phi(i+ct)
\end{equation}
with 
$$\phi'\geq 0\quad\mbox{ and }\quad \phi\ \mbox{ is bounded}.$$
Such a solution, $\phi,$ is called a traveling wave solution and $c$ denotes its
velocity of propagation. 
From (\ref{e1}) and (\ref{e3}), it is equivalent to look for solutions $\phi$ of
\begin{equation}\label{e6}
c\phi'(z)=\phi(z+1)-2\phi(z)+\phi(z-1)+f(\phi(z))+\sigma
\end{equation}
with $z=i+ct.$ 
For such a model, and under certain conditions on $f,$ we show the existence of traveling waves for each value of $\sigma$ in an interval $[\sigma^{-},\sigma^{+}]$ (see Theorem \ref{t2}). We distinguish three "phases": $\sigma\in (\sigma^-,\sigma^+)$ for bistable nonlinearities, $\sigma=\sigma^+$ for positive monostable nonlinearity and $\sigma=\sigma^-$ for negative monostable nonlinearity. Those three phases match together in a unified picture that we call the velocity  diagram. On Figure \ref{f3}, the diagram shows the nondecreasing velocity function 
$c=c(\sigma)$ with respect to the driving force $\sigma\in (\sigma^-,\sigma^+)$, with vertical branches $\displaystyle c\ge c^+$ for $\sigma=\sigma^{-}$, and $\displaystyle c\le c^-$ for $\sigma=\sigma^{+}$. Here the critical velocities $c^+,c^-$ are limits in the monostable case of the velocities in the bistable cases
$$c^+:=\lim_{\sigma^+>\sigma \to \sigma^+}c(\sigma),\quad c^-:=\lim_{\sigma^-<\sigma \to \sigma^-}c(\sigma)$$

The goal of this paper is to show that such behaviour arises  in a  framework which is much more general than \eqref{e6}: the fully nonlinear framework.
To this end, given a real function $F$ (whose properties will be specified later),
we consider the following generalized equation with parameter $\sigma\in\R$
\begin{equation}\label{e8}
c\phi'(z)=F(\phi(z+r_{0}),\phi(z+r_{1}),...,\phi(z+r_{N}))+\sigma,
\end{equation}
where $N\geq 0$ and $r_{i}\in\R$ for $i=0,...,N$ such that
\begin{equation}\label{e9}
r_{0}=0\quad\mbox{and}\quad r_{i}\neq r_{j}\ \mbox{ if }\ i\neq j,
\end{equation}
which does not restrict the generality. In \eqref{e8}, we are looking for both the profile $\phi$ and the velocity $c.$\\

Traveling waves were studied also for discrete bistable reaction-diffusion equations (see for instance \cite{CCHM,CGW}). 
See also \cite{AFM,HMSV} and the references therein. In the monostable case, we distinguish \cite{HZ} (for nonlocal non-linearities with integer shifts) and \cite{CDM,LWL,We1,Y} (for problems with linear nonlocal part and with integer shifts also). See also \cite{GW08} for particular monostable nonlinearities with irrational shifts. We also refer to \cite{GH,CFG,GW10,CG02,CG03,HN3,ZHH} for different positive monostable nonlinearities. In the monostable case, we have to underline the work of Hudson and Zinner \cite{HZ} (see also \cite{ZHH}), where they proved the existence of a branch of solutions $c\geq c^{*}$ for general Lipschitz nonlinearities (with possibly an infinite number of neighbors $N=+\infty$, and possibly $p$ types of different particles, while $p=1$ in our study) but with integer shifts $r_i\in \Z.$ However, they do not state the nonexistence of solutions for $c<c^{*}.$ Their method of proof relies on an approximation of the equation on a bounded domain (applying Brouwer's fixed point theorem) and an homotopy argument starting from a known solution. The full result is then obtained as the size of the domain goes to infinity.
Here we underline that our results hold for the fully nonlinear case with real shifts $r_{i}\in\R.$

Several approaches were used to construct traveling waves for discrete monostable dynamics. We already described the homotopy method of Hudson and Zinner \cite{HZ}. 
In a second approach, Chen and Guo \cite{CG02} proved the existence of a solution starting from an approximated problem. They constructed a 
fixed point solution of an integral reformulation (approximated on a bounded domain) using the monotone iteration method (with sub and supersolutions). 
This approach was also used to get the existence of a solution in \cite{FGS,CG03,GW08,GW10}. A third approach based on recursive method for monotone discrete in time dynamical systems was used by Weinberger et al. \cite{LWL,We1}. See also \cite{Y}, where this method 
is used to solve problems with a linear nonlocal part. 
In a fourth approach \cite{GH}, Guo and Hamel used global space-time sub and supersolutions to prove the existence of a solution for periodic monostable equations.  

There is also a wide literature about the uniqueness and the asymptotics at infinity of a solution for a monostable non-linearities, see for instance \cite{CFG,HLM} (for a degenerate case), \cite{CG02,CG03} and the references therein. 
Let us also mention that certain delayed reaction diffusion equations with some KPP-Fisher non-linearities do not admit traveling waves (see for example \cite{FGS,ZHH}).

Finally, we mention that our method opens new possibilities to be adapted to more general problems. 
For example, we can think to adapt our 
approach to a case with possibly $p$ types of different particles similar to \cite{FIM12}. The case with an infinite number of neighbors $N=+\infty$ could be also studied. We can also think to study fully nonlinear parabolic equations.\\

The present work has been already announced in a preprint \cite{AM} that was accessible since 2014 and also in the PhD thesis in 2014 of the first author.
Unfortunately, the life conditions of the two authors did not permit the submission  to publication of the manuscript. 
The present paper corresponds to parts I and II of \cite{AM}. Part III of \cite{AM} will be presented in another work \cite{AM-II}, where we study the general Lipschitz monostable cases which can no longer be seen as a limit case of the bistable case.

Notice also that part of our work can be seen as an extension to Lipschitz discrete dynamics of some results that hold true for classical reaction diffusion equations (see \cite{AM-cras}).\\

\subsection{Main results}\label{sbbs}

In this subsection, we consider equation \eqref{e8} with a constant parameter $\sigma\in\R$ and $F:\R^{N+1}\to\R.$ We are interested in the velocities $c$ associated to $\sigma$ (that we call roughly speaking the ``velocity function").

For $\sigma$ belonging to some interval $[\sigma^{-},\sigma^{+}],$ we prove the existence of a traveling wave and we study the variation of its velocity $c$ with respect to $\sigma.$

\noindent
Let $E=(1,...,1),\ \Theta=(\theta,...,\theta)\in\R^{N+1}$ with $\theta \in (0,1)$ and assume that the function $F$ satisfies:\\

\noindent{\bf Assumption $\boldsymbol{(\tilde{A}_{C^{1}})}$:}
\begin{description}
\item[\quad{\bf Regularity:}] $F$ is globally Lipschitz continuous over $\R^{N+1}$ and $C^{1}$ over a neighborhood in $\R^{N+1}$ of the two intervals $]0,\Theta[$ and $]\Theta,E[.$
\item[\quad{\bf Monotonicity:}] $F(X_{0},...,X_{N})$ is non-decreasing w.r.t. each $X_{i}$ for $i\neq 0.$
\item[\quad{\bf Periodicity:}] $F(X_{0}+1,...,X_{N}+1)=F(X_{0},...,X_{N})$ for every $X=(X_{0},...,X_{N})\in\R^{N+1}.$
\end{description}

\noindent
Notice that, since $F$ is periodic in $E$ direction, then $F$ is $C^{1}$ over a neighborhood of $\R E\backslash(\Z E\ \cup\ \Z\Theta).$ 

\medskip

\noindent{\bf Assumption $\boldsymbol{(\tilde{B}_{C^{1}})}$:}\\
Define $f(v)=F(v,...,v)$ such that: 
\begin{description}
\item[\quad{\bf Bistability:}] $f(0)=f(1)$ and there exists $\theta\in(0,1)$ such that 
$$\left\{
\begin{aligned}
&f'>0\quad\mbox{on}\quad(0,\theta)\\
&f'<0\quad\mbox{on}\quad(\theta,1)\\
\end{aligned}
\right.$$
\end{description}

\begin{figure}[!ht]
\centering\epsfig{figure=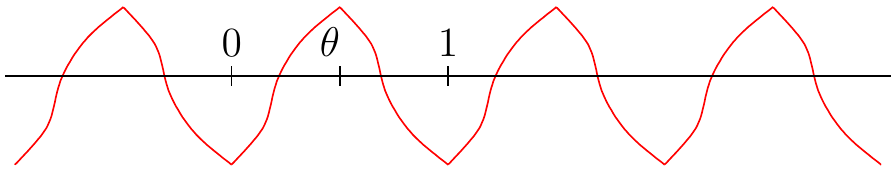,width=80mm}
\caption{Bistable nonlinearity $f$\label{f2}}
\end{figure}

\noindent
See Figure \ref{f2} for an example of $f$ satisfying $(\tilde{B}_{C^{1}}).$ Notice that assumptions $(\tilde{A}_{C^{1}})$ and $({\tilde B}_{C^{1}})$ holds true in particular for the Frenkel-Kontorova model (\ref{e1}). 

\begin{theo}\label{t2}{\bf (Velocity diagram for  traveling waves)}\\
Under assumptions $(\tilde{A}_{C^{1}})$ and $({\tilde B}_{C^{1}}),$ 
define $\sigma^{\pm}$ as 
$$\left\{
\begin{aligned}
&\sigma^{+}=-\min f\\
&\sigma^{-}=-\max f.
\end{aligned}
\right.$$
Associate for each $\sigma\in[\sigma^{-},\sigma^{+}]$ the (unique) solutions $m_{\sigma}\in[\theta-1,0]$ and $b_{\sigma}\in[0,\theta]$ of $f(s)+\sigma=0.$ Then consider the following equation
\begin{equation}\label{e15}
\left\{
\begin{aligned}
&c\phi'(z)=F(\phi(z+r_{0}),\phi(z+r_{1}),...,\phi(z+r_{N}))+\sigma\quad\mbox{on}\quad\R\\
&\phi\mbox{ is non-decreasing over }\R\\
&\phi(-\infty)=m_{\sigma}\quad\mbox{and}\quad\phi(+\infty)=m_{\sigma}+1,
\end{aligned}
\right.
\end{equation}
{\bf $1$- Bistable case: traveling waves for $\sigma\in(\sigma^{-},\sigma^{+})$}\\
We have
\begin{itemize}
\item[\quad\quad{\bf $(i)$}]{\bf (Existence of a traveling wave)}\\
For any $\sigma\in(\sigma^{-},\sigma^{+}),$ there exists a unique real $c:=c(\sigma),$ such that there exists a function $\phi_{\sigma}:\R\to\R$ solution of (\ref{e15}) in the viscosity sense. 
\item[\quad\quad{\bf $(ii)$}]{\bf (Continuity and monotonicity of the velocity function)}\\
The map $$\sigma\mapsto c(\sigma)$$ is
continuous on $(\sigma^{-},\sigma^{+})$ and there exists a constant $K>0$ such that the function $c(\sigma)$ is non-decreasing and satisfies $$\frac{dc}{d\sigma}\geq K|c| \quad\mbox{on }\ (\sigma^{-},\sigma^{+})$$ in the viscosity sense. 
In addition, there exist real numbers $c^{-}\leq c^{+}$ such that 
$$\lim_{\sigma^-<\sigma\to\sigma^{-}}c(\sigma)=c^{-}\quad\mbox{and}\quad\lim_{\sigma^+>\sigma\to\sigma^{+}}c(\sigma)=c^{+}.$$
Moreover, either $c^{-}=0=c^{+}$ or $c^{-}<c^{+}.$
\end{itemize}
{\bf$2$- Monostable cases: vertical branches for $\sigma=\sigma^{\pm}$}\\ 
We have
\begin{itemize}
\item[\quad\quad{\bf $(i)$}]{\bf (Existence of traveling waves for $c\geq c^{+}$ when $\sigma=\sigma^{+}$)}\\ Let $\sigma=\sigma^{+},$ then for every
$c\geq c^{+}$ there exists a traveling wave $\phi$ solution of
\begin{equation}\label{e17}
\left\{
\begin{aligned}
&c\phi'(z)=F(\phi(z+r_{0}),\phi(z+r_{1}),...,\phi(z+r_{N}))+\sigma^{+}\quad\mbox{on}\quad\R\\
&\phi\mbox{ is non-decreasing over }\R\\
&\phi(-\infty)=0=m_{\sigma^{+}}\quad\mbox{and}\quad\phi(+\infty)=1.
\end{aligned}
\right.
\end{equation}
Moreover, for any $c<c^{+},$ there is no solution $\phi$ of (\ref{e17}).  
\item[\quad\quad{\bf $(ii)$}]{\bf (Existence of traveling waves for $c\leq c^{-}$ when $\sigma=\sigma^{-}$)}\\ Let $\sigma=\sigma^{-},$ then for every $c\leq c^{-},$ there exists a traveling wave $\phi$ solution of
\begin{equation}\label{e17'}
\left\{
\begin{aligned}
&c\phi'(z)=F(\phi(z+r_{0}),\phi(z+r_{1}),...,\phi(z+r_{N}))+\sigma^{-}\quad\mbox{on}\quad\R\\
&\phi\mbox{ is non-decreasing over }\R\\
&\phi(-\infty)=\theta-1=m_{\sigma^{-}}\quad\mbox{and}\quad\phi(+\infty)=\theta.
\end{aligned}
\right.
\end{equation}
Moreover, for any $c>c^{-},$ there is no solution $\phi$ of (\ref{e17'}).
\end{itemize}
\end{theo}

We have to mention that Theorem \ref{t2}-$1$ $(i)$ is already proved in \cite{AFM} (see \cite[Proposition 2.3]{AFM}). 
Hence our contribution consists in the remaining parts of theorem. The originality of our work is probably more in the statement of the theorem than in the proof itself.\\

\noindent
\begin{figure}[!ht]
\centering\epsfig{figure=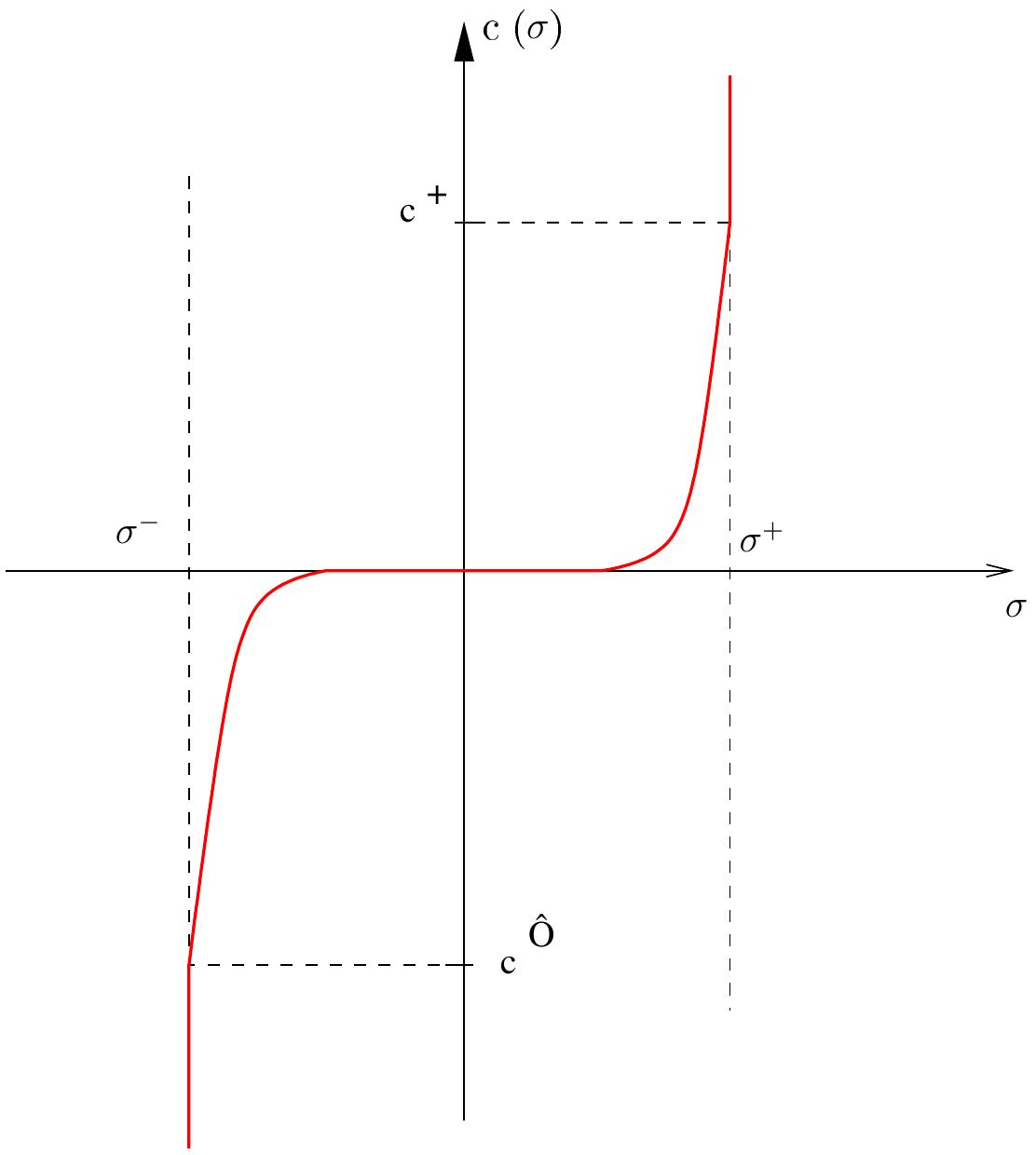,width=80mm}
\caption{Schematic velocity diagram: the velocity function $c(\sigma)$ with vertical branches at $\sigma=\sigma^{\pm}.$\label{f3}}
\end{figure}

Notice that there are no monotone traveling waves solutions for $\sigma\not\in [\sigma_-,\sigma_+]$ because there are no solutions of $f(\cdot)+\sigma=0$ in that case.

Notice also that from Proposition \ref{pdvs} in the Appendix, we know that for the Frenkel-Kontorova model \eqref{e1}, we have $\sigma^{\pm}=\pm\beta$ and $c^{+}>0>c^{-}$, and Figure \ref{f3} illustrates the velocity diagram. Moreover for this particular model the velocity function $c(\sigma)$  has a plateau at the level $c=0$ in particular if $|\sigma|<\beta-1$. Another example where the velocity is $c=0$ is in the special case where $F(X_0,\dots,X_N)=f(X_0)$ and $\sigma\in (\sigma^-,\sigma^+)$. This last example  is also a situation where the traveling waves are discontinuous, constant on some half lines, and not unique (even up to translations), while we have uniqueness of the velocity $c=0$ for $\sigma\in (\sigma^-,\sigma^+)$. This example also shows that the strong maximum principle may not hold in our general setting.

Notice also that when $c\not=0$, then the solutions $\phi$ are classical, i.e. satisfy $\phi \in C^{1,1}$. On the contrary, when $c=0$, the profile $\phi$ may be discontinuous, as we have seen. Indeed, the notion of viscosity solution can be used for any value of $c\in \R$. It has also be seen in \cite{AFM}, that for monotone profiles $\phi$, there is equivalence to be a viscosity solution and to be a solution almost everywhere.\\

\noindent In view of Theorem \ref{t2}, we can ask the following (see also \cite{CCHM}):\\

\noindent{\bf Open question}:\\
For a general $F,$ what is the precise behavior of the function $c(\sigma)$ close to the boundary of the plateau $c=0$ (when it does exist) and close to $\sigma^{+}$ and $\sigma^{-}?$\\

As a notation, we set for a general function $h$:
$$F((h(z+r_{i}))_{i=0,...,N})=F(h(z+r_{0}),h(z+r_{1}),...,h(z+r_{N}))$$ and we define
\begin{equation}\label{e2}
r^{*}=\max_{i=0,...,N}|r_{i}|.
\end{equation}
In the rest of the paper, we will use the notation introduced in Theorem \ref{t2}.

\subsection{Organization of the paper}

In Section \ref{S2}, we give/recall some  results concerning viscosity solutions and stability of solutions that will be useful in the next sections.
In Section \ref{S3}, we show that for large velocities $c>>1$, we have existence of traveling waves for the parameter $\sigma=\sigma^+$ (and also give a similar result for $\sigma=\sigma^-$).
In Section \ref{S4}, we study the properties of the velocity function $\sigma\mapsto c(\sigma)$ for the range $\sigma\in (\sigma^-,\sigma^+)$. We also use the monotonicity of the velocity function $c(\sigma)$ to define the limits $c^\pm$ of the velocities for the values $\sigma=\sigma^\pm$.
In Section \ref{S5}, we show that there are no traveling waves for $c<c^+$ and $\sigma=\sigma^+$, and construct a branch of traveling waves for each velocity $c\ge c^+$ for $\sigma=\sigma^+$. We get a similar result for $\sigma=\sigma^-$. Finally at the end of Section \ref{S5}, we give the proof of Theorem \ref{t2}, as a corollary of all the previous results.

\section{Useful results}\label{S2}

We recall here some useful results. Some of these results are contained in \cite{AFM}.\\

We first recall the  notion of viscosity solutions that we use in this work. To this end,
we recall that the upper and lower semi-continuous envelopes, $u^{*}$ and $u_{*},$ of a locally bounded function $u$ are defined as
$$u^{*}(y)=\limsup_{x\rightarrow y}u(x)\quad\mbox{and}\quad u_{*}(y)=\liminf_{x\rightarrow y}u(x).$$
and that $u$ is upper semi-continuous if and only if $u=u^*$ (and similarly $u$ is lower semi-continuous if and only if $u=u_*$).

\begin{defi}\label{d1}{\bf (Viscosity solution)}\\
Let $I=I'=\R$ (or $I=(-r^{*},+\infty)$ and $I'=(0,+\infty)$) and $u:I\rightarrow\R$ be a locally bounded function, $c\in\R$ and $F$ defined on $\R^{N+1}$.
\begin{itemize}
\item[-]The function $u$ is a subsolution (resp. a supersolution)
  on $I'$ of 
\begin{equation}\label{e47}
cu'(x)=F((u(x+r_{i}))_{i=0,...,N})+\sigma,
\end{equation}
 if $u$ is upper semi-continuous
  (resp. lower semi-continuous) and if for all
  test function $\psi\in C^{1}(I)$ such that $u-\psi$ attains a local
  maximum (resp. a local minimum) at $x^{*}\in I',$ we
  have $$c\psi'(x^{*})\leq
  F((u(x^{*}+r_{i}))_{i=0,...,N})+\sigma\mbox{\quad\Big(resp. }c\psi'(x^{*})\geq F((u(x^{*}+r_{i}))_{i=0,...,N})+\sigma\Big).$$
\item[-]A function $u$ is a viscosity solution of (\ref{e47}) on $I'$ if $u^{*}$ is a subsolution and $u_{*}$ is a supersolution on $I'.$
\end{itemize}
\end{defi}

Next, we state Perron's method  to construct solutions.

\begin{pro}\label{p6}{\bf(Perron's method (\cite[Proposition $2.8$]{FIM}))}\\
Let $I=(-r^{*},+\infty)$ and $I'=(0,+\infty)$ and $F$ be a function satisfying $(\tilde{A}_{C^1})$. Let $u$ and $v$ defined on $I$ satisfying $$u\leq v\quad\mbox{on}\quad I,$$ such that $u$ and $v$ are respectively a sub and a supersolution of (\ref{e47}) on $I'.$ Let $\mathcal{L}$ be the set of all functions $\tilde{v}:I\to\R,$ such that $u\leq\tilde{v}$ over $I$ with $\tilde{v}$ supersolution of (\ref{e47}) on $I'.$ For every $z\in I,$ let $$w(z)=\inf\{\tilde{v}(z)\quad\mbox{such that}\quad\tilde{v}\in\mathcal{L}\}.$$ Then $w$ is a viscosity solution of (\ref{e47}) over $I'$ satisfying $u\leq w\leq v$ over $I.$
\end{pro}

The following result is important and meaningful in our work.

\begin{lem}\label{l8}{\bf(Equivalence between viscosity and a.e. solutions, \cite[Lemma $2.11$]{AFM})}\\
Let $F$ satisfying assumption $(\tilde{A}_{C^1})$. Let $\displaystyle\phi:\R\rightarrow\R$ be a non-decreasing function. Then $\phi$ is a viscosity solution of
$$c\phi'(x) =F((\phi(x+r_{i}))_{i=0,...,N})+\sigma\quad\mbox{on}\quad\R,$$
if and only if $\phi$ is an almost everywhere solution of the same equation.
\end{lem}

Having this result in hands, we have the following useful criterion to pass to the limit.

\begin{pro}\label{pro::r2}{\bf (Stability by passage to the limit)}\\
Let $F$ satisfying assumption $(\tilde{A}_{C^1})$. Given $a<b$, let $\displaystyle\phi_n:I:=(a-r^*,b+r^*) \rightarrow\R$ be a non-decreasing  viscosity solution of
$$c_n\phi_n'(x) =F((\phi_n(x+r_{i}))_{i=0,...,N})+\sigma\quad\mbox{on}\quad I':=(a,b)$$
satisfying the bounds
$$|\phi_n|_{L^\infty(I)}\le C,\quad |c_n|\le C$$
Then up to a subsequence, we have
$$\phi_n \to \phi \quad \mbox{a.e. on $I$},\quad c_n \to c$$
and $\phi$ is a viscosity solution of
$$c\phi'(x) =F((\phi(x+r_{i}))_{i=0,...,N})+\sigma\quad\mbox{on}\quad I'$$
\end{pro}

\noindent {\bf Proof of Proposition \ref{pro::r2}}\\
The existence of a subsequence converging almost everywhere follows from classical Helly's theorem for monotone functions.
The remaining part of the argument follows from the equivalence between viscosity solutions and almost everywhere solutions when $c=0$.
In the case $c\not=0$, we get bounds on $|\phi_n|_{C^1(I)}\le C'$, and the result follows for instance from the classical stability of viscosity solutions
(or also by a direct argument for ODEs).\\

\section{Vertical branches for large velocities}\label{S3} 

In this section, we want to build traveling waves for large velocities in the case $\sigma=\sigma^\pm$.
To this end, we focus on the case $\sigma=\sigma^+$ (and we will see later that the case $\sigma_-$ is similar).
To simplify the presentation, we assume (without loss of generality) that
$$\sigma^+=0$$
which means under assumptions $(\tilde A_{C^1})$ and $(\tilde B_{C^1})$ that
$$f(0)=f(1)=0,\quad f>0 \quad \mbox{on}\quad (0,1)$$

For $\sigma=\sigma^+=0$, we now want to study solutions $(c,\phi)$  of the equation
\begin{equation}\label{e11}
\left\{
\begin{aligned}
&c\phi'(z)=F(\phi(z+r_{0}),\phi(z+r_{1}),...,\phi(z+r_{N}))\quad\mbox{on}\quad\R\\
&\phi\mbox{ is non-decreasing over }\R\\
&\phi(-\infty)=0\quad\mbox{and}\quad\phi(+\infty)=1.
\end{aligned}
\right.
\end{equation}
and show the existence of solutions for $c>>1$ large enough.\\

To this end, we first need the following definition

\begin{defi}\label{supsol-equa}{\bf(Supersolution of \eqref{e11})}\\
We say that $(c,\psi)$ is a supersolution of \eqref{e11} if $(c,\psi)$ satisfies
$$
\left\{
\begin{aligned}
&c\psi'(z)\geq F((\psi(z+r_{i}))_{i=0,...,N})\quad\mbox{on}\quad\R\\
&\psi\mbox{ is non-decreasing over }\R\\
&\psi(-\infty)=0\quad\mbox{and}\quad\psi(+\infty)=1.
\end{aligned}
\right.
$$
\end{defi}

Given a general supersolution $\psi$, it may be difficult to build a solution of (\ref{e11}) if $\psi$ vanishes on a left half line $(-\infty,0)$, which may happen when we do not have strong maximum principle (which is not assumed to hold for general $F$ as considered here).
On the contrary, when $\psi$ is positive, we have the following useful result.

\begin{pro}\label{sup->sol}{\bf(Solution of \eqref{e11} if it admits a positive supersolution)}\\
Consider a function $F$ satisfying $(\tilde A_{C^1})$ and $(\tilde B_{C^1})$ and assume that $\sigma_+=0$. Assume that there exists a supersolution $(c,\psi)$ of \eqref{e11} with  $\psi>0$. Then there exists a traveling wave $\phi$ such that $(c,\phi)$ is a solution of \eqref{e11}.  
\end{pro}

\noindent {\bf Proof of Proposition \ref{sup->sol}}\\
Up to translation, we can assume that $(0,1)\ni \theta \in [\psi_*(0),\psi^*(0)]$. We will construct a solution using Perron's method.\\
\noindent{\bf Step $1$: construction of a subsolution}\\
Consider the constant function $\overline{\psi}=\varepsilon$ with $\varepsilon>0$ small enough fixed. Then 
$$0=c\overline{\psi}'(x)\leq F((\overline{\psi}(x+r_{i}))_{i=0,...,N})=f(\varepsilon).$$ Hence $(c,\overline{\psi})$ is a subsolution of
\begin{equation}\label{ode-F}
cw'(x)=F((w(x+r_{i}))_{i=0,...,N})\quad\mbox{on}\quad\R.
\end{equation}
\noindent{\bf Step $2$: construction of solution on a half line}\\
Since $\psi(-\infty)=0,$ $\psi(0)=\theta>0$ and $\psi$ is non-decreasing and positive,
then for $\varepsilon \in (0,\theta)$ fixed, we can define $k_{\varepsilon}<0$ such that 
$$\varepsilon\in [\psi_*(k_\varepsilon),\psi^*(k_\varepsilon)]\quad\mbox{and}\quad\psi>\varepsilon\quad\mbox{on}\quad (k_{\varepsilon},+\infty).$$
Then using Perron's method (Proposition \ref{p6}), there exists a solution $\phi_{\varepsilon}$ of \eqref{ode-F} on $(r^{*}+k_{\varepsilon},+\infty)$ such that
$$\varepsilon\leq\phi_{\varepsilon}\leq\psi\quad\mbox{on}\quad(k_{\varepsilon},+\infty).$$
\noindent{\bf Step $3$: $\phi_{\varepsilon}$ is non-decreasing on $\left(k_{\varepsilon},+\infty\right).$}\\
Define for $x\in(k_{\varepsilon},+\infty)$ the function $$\overline{\phi}(x):=\inf_{p\geq 0}\phi_{\varepsilon}(x+p).$$ Clearly, since $\varepsilon\leq \phi_{\varepsilon}(x+p)$ for all $p\geq0$ and $x\in\left(k_{\varepsilon},+\infty\right),$ we get $\varepsilon\leq\overline{\phi}(x)\leq \phi_{\varepsilon}(x)\leq\psi(x)$ for all $x\in\left(k_{\varepsilon},+\infty\right).$ On the other hand, for all $p\geq 0,$ $\phi_{\varepsilon}(x+p)$ is a solution of \eqref{ode-F} over $\left(r^{*}+k_{\varepsilon},+\infty\right),$
then $(\overline{\phi})_{*}$ is supersolution of \eqref{ode-F} over $\left(k_{\varepsilon}+r^{*},+\infty\right)$ (as it is classical for viscosity solutions). Moreover, we have $\varepsilon\leq(\overline{\phi})_{*}\leq\psi$ on $(k_\varepsilon,+\infty)$.\\
Now because $\phi_{\varepsilon}$ itself is defined in Perron's method as an infimum over all supersolutions above $\varepsilon$, we deduce that
$$\phi_{\varepsilon}\le \overline{\phi}$$
while we have the reverse inequality by definition of $\overline{\phi}$. Hence we get
$$\phi_{\varepsilon}= \overline{\phi}$$
which shows that $\phi_{\varepsilon}$ is non-decreasing over $\left(k_{\varepsilon},+\infty\right)$.\\
\noindent{\bf Step $4$: passing to the limit $\varepsilon\to0$}\\
Since $\phi_{\varepsilon}$ is a non-decreasing solution  on $(r^{*}+k_{\varepsilon},+\infty),$
then $\phi_{\varepsilon}(+\infty)$ has to solve $f(x)=0$. Moreover we have $0<\varepsilon\leq\phi_{\varepsilon}\leq\psi\leq1$ over $(k_{\varepsilon},+\infty)$, and we conclude that $$\phi_{\varepsilon}(+\infty)=1.$$
Now from the fact that $\theta\in [\psi_*(0),\psi^*(0)]$, we deduce that there exists $x_\varepsilon\ge 0$ such that
$$\theta \in [(\phi_\varepsilon)_*(x_\varepsilon),(\phi_\varepsilon)^*(x_\varepsilon)]$$
Moreover, because $\psi>0$, we deduce that
$$k_\varepsilon\to -\infty\quad \mbox{as}\quad \varepsilon\to 0^+$$
Then we can shift $\phi_\varepsilon$, defining
$$\tilde \phi_\varepsilon(x):= \phi_\varepsilon(x+x_\varepsilon)$$
which satisfies
$$\theta\in [(\tilde \phi_\varepsilon)_*(0),(\tilde \phi_\varepsilon)^*(0)]$$
and $\tilde \phi_\varepsilon$ is a solution of the equation on $(r^*+\tilde k_\varepsilon,+\infty)$ with
$$\tilde k_\varepsilon := k_\varepsilon-x_\varepsilon \to -\infty \quad \mbox{as}\quad \varepsilon\to 0^+$$
Using the stability of solutions (Proposition \ref{pro::r2}), we deduce that we can pass to the limit $\tilde \phi$ of $\tilde \phi_\varepsilon$ (with at least convergence almost everywhere), and that the limit is still a solution, i.e. satisfies
$$\left\{
\begin{aligned}
&c\tilde \phi'(z)=F(\tilde\phi(z+r_{0}),\tilde\phi(z+r_{1}),...,\tilde\phi(z+r_{N}))\quad\mbox{on}\quad\R\\
&\tilde\phi\mbox{ is non-decreasing over }\R\\
&\theta\in [(\tilde \phi)_*(0),(\tilde \phi)^*(0)]\\
& 0\le \tilde \phi \le 1
\end{aligned}
\right.$$
which implies that
$$\tilde \phi(-\infty)=0,\quad \tilde \phi(+\infty)=1$$
This shows that $\tilde \phi$ is a solution of (\ref{e11}) and ends the proof of the proposition.\\

\begin{cor}\label{sol->sup}{\bf(Half line of solutions)}\\
Consider a function $F$ satisfying $(\tilde A_{C^1})$ and $(\tilde B_{C^1})$ and assume that $\sigma^+=0$. Assume that there exists a supersolution $(c,\phi)$ of \eqref{e11} with  $\phi>0$. Then for all $\tilde{c}\geq c$ there exists a solution $\tilde{\phi}$ of \eqref{e11}.
\end{cor}

\noindent {\bf Proof of Corollary \ref{sol->sup}}\\
We simply notice that for $\tilde c\ge c$, from the monotonicity of $\phi$, we deduce in the sense of viscosity
$$\tilde{c}\phi'(z)\geq c\phi'(z)=F((\phi(z+r_{i}))_{i=0,...,N}).$$
Because $\phi$ is a positive supersolution for the velocity $\tilde c$, we deduce from Proposition \ref{sup->sol} that there exists a solution $\tilde \phi$ for the velocity $\tilde c$. This ends the proof of the corollary.\\

\begin{pro}\label{p0}{\bf(Existence of traveling waves for $c>>1$)}\\
Consider a function $F$ satisfying $(\tilde A_{C^1})$ and $(\tilde B_{C^1})$ and assume that $\sigma^+=0$. 
Then there exists some $c_L>0$ large enough such that for any $c\ge c_L$, there exists some  traveling wave $\phi$ solution of (\ref{e11}) with velocity $c$.
\end{pro}

\noindent{\bf Proof of Proposition \ref{p0}}\\
\noindent{\bf Step 1: preliminary}\\
In order to build a solution $\phi$ for large velocity $c$ of
\begin{equation}\label{e19}
c\phi'(x)=F(\phi(x+r_{0}),\phi(x+r_{1}),...,\phi(x+r_{N}))\quad\mbox{on}\quad\R,
\end{equation}
it is convenient to set
$$h(z):=\phi(cz)\quad \mbox{and}\quad \varepsilon:=\frac{1}{c}$$
where we now look for some $h$ solution of
\begin{equation}\label{e20}
h'(z)=F\left(\left(h\left(z+\varepsilon r_{i}\right)\right)_{i=0,...,N}\right)\ \mbox{ on }\ \R.
\end{equation}
Then it is natural to look for the unique solution of the associated ODE
$$h'_{0}=F(h_{0},...,h_{0})=f(h_{0})\geq 0, \quad h_0(0)=\theta$$
which satisfies
$$h_0(-\infty)=0,\quad h_0(+\infty)=1$$
We now consider
$$h_\varepsilon(z):=h_0(a_\varepsilon z)\quad \mbox{with}\quad a_\varepsilon:=1+M\varepsilon$$
and want to show below that for $M>0$ large enough, and $\varepsilon$ small enough, the function $h_\varepsilon$ is a supersolution of (\ref{e20}).\\
\noindent{\bf Step $2$: estimate on $h_0\left(z+b\right)$}\\
Recall that $f \in Lip([0,1])$ and then that $h_0\in W^{2,\infty}(\R)$. We have
$$h_0(z+b)=h_0(z)+b \int_0^1 h_0'(z+bt)\ dt$$
and (at least almost everywhere, and indeed everywhere except at $z=0$ where $h_0(0)=\theta$)
$$h_0''=f'(h_0)h_0'$$
ie
$$(\ln h_0')' = f'(h_0)$$
which implies
$$h_0'(z+bt) \le h_0'(z)  e^{|f'|_{L^\infty}|b|t}$$
and then
$$h_0(z+b)\le h_0(z)+ h_0'(z) \cdot g(b)\quad \mbox{with}\quad g(b):=|b| e^{|f'|_{L^\infty}|b|}$$
\noindent{\bf Step $3$: checking that $h_\varepsilon$ is a supersolution for $M>0$ large enough}\\
We deduce that
$$h_\varepsilon(z+\varepsilon r_i) =h_0(a_\varepsilon (z+ \varepsilon r_i))\le h_\varepsilon(z)+ \frac{1}{a_\varepsilon}\cdot h_\varepsilon'(z) g(\varepsilon a_\varepsilon r_i) $$
Because of the monotonicities of $F$ and the fact that $F$ is Lipschitz continuous, we deduce that there exists a constant $K>0$ such that
$$F\left(\left(h_\varepsilon\left(z+\varepsilon r_{i}\right)\right)_{i=0,...,N}\right) - F\left(\left(h_\varepsilon\left(z\right)\right)_{i=0,...,N}\right)\le K \max_{i=0,\dots,N} \left\{\frac{1}{a_\varepsilon}\cdot h_\varepsilon'(z) g(\varepsilon a_\varepsilon r_i)\right\}$$
Now
$$F\left(\left(h_\varepsilon\left(z\right)\right)_{i=0,...,N}\right)=f(h_\varepsilon(z))=h_0'(a_\varepsilon z)=\frac{1}{a_\varepsilon} h_\varepsilon'(z)$$
which shows that
$$\mu_\varepsilon \cdot h_\varepsilon'(z) \ge F\left(\left(h_\varepsilon\left(z+\varepsilon r_{i}\right)\right)_{i=0,...,N}\right)\quad \mbox{with}\quad \mu_\varepsilon:=\frac{1+K \max_{i=0,\dots,N}g(\varepsilon a_\varepsilon r_i)}{a_\varepsilon} $$
We now claim that 
$$\mu_\varepsilon \le 1$$
for $\varepsilon$ small enough and $M>0$ large enough. Indeed, for $\varepsilon$ small enough, we can insure
$$a_\varepsilon \le 2,\quad \varepsilon\le 1$$
and then
$$g(\varepsilon a_\varepsilon r_i) \le \varepsilon (2r^*) e^{2r^* |f'|_{L^\infty}}$$
which gives
$$\mu_\varepsilon \le \frac{1+ \varepsilon K (2r^*) e^{2r^* |f'|_{L^\infty}}}{1+\varepsilon M}\le 1\quad \mbox{for}\quad M\ge K(2r^*) e^{2r^* |f'|_{L^\infty}}$$
This finally implies that
$$h_\varepsilon'(z) \ge F\left(\left(h_\varepsilon\left(z+\varepsilon r_{i}\right)\right)_{i=0,...,N}\right)$$
ie that $h_\varepsilon$ is a supersolution.\\
\noindent{\bf Step $4$: Conclusion}\\
We now see that $\phi_\varepsilon$ defined by
$$h_\varepsilon(z)=:\phi_\varepsilon(cz)\quad \mbox{with}\quad c=\frac{1}{\varepsilon}$$
is a positive supersolution of the original equation (\ref{e19}) for the velocity $c$. Hence we can apply Proposition \ref{e11}
which shows the existence of a traveling wave for the velocity $c$.
This ends the proof of the proposition.\\

\begin{cor}\label{cor::r11}{\bf (Existence of traveling waves for large negative $c$)}\\
Consider a function $F$ satisfying $(\tilde A_{C^1})$ and $(\tilde B_{C^1})$ and assume that $\sigma=\sigma^-$. 
Then there exists some $c_L'>0$ large enough such that for any $c\le -c_L'$, there exists some  traveling wave $\phi$ solution of (\ref{e17'}) with velocity $c$.
\end{cor}

\noindent {\bf Proof of Corollary \ref{cor::r11}}\\
The result follows from the fact that $\phi$ is a solution of 
$$c\phi'(z)=F((\phi(z+r_{i}))_{i=0,...,N})+\sigma^{-}\quad\mbox{over}\quad\R$$
if and only if 
$$\bar{\phi}(z)=1-\phi(-z)$$ is a solution of 
$$\bar c \bar \phi'(z)=\bar F((\bar \phi(z+\bar r_i))_{i=0,\dots,N})+ \bar \sigma^+$$
with
$$\left\{
\begin{aligned}
&\bar{F}(X_{0},...,X_{N})=-F((1-X_{i})_{i=0,...,N})\\
&\bar{c}=-c\\
&\bar{r}_{i}=-r_{i}\\
&\bar{\sigma}^{+}=-\sigma^{-}\\
\end{aligned}
\right.$$
This implies in particular that
$$\bar f(v):=\bar F(v,\dots,v)\ge 0$$
and $\bar F$ still satisfies conditions $(\tilde A_{C^1})$ and $(\tilde B_{C^1})$ with $\theta$ replaced by $\bar \theta:=1-\theta$. Hence we can apply Proposition \ref{p0} which leads to the result.
This ends the proof of the corollary.\\

\section{Traveling waves for $\sigma\in (\sigma^-,\sigma^+)$ and properties of the velocity}\label{S4} 

In this section, we consider equation (\ref{e15}), namely
\begin{equation}\label{e15bis}
\left\{
\begin{aligned}
&c\phi'(z)=F(\phi(z+r_{0}),\phi(z+r_{1}),...,\phi(z+r_{N}))+\sigma\quad\mbox{on}\quad\R\\
&\phi\mbox{ is non-decreasing over }\R\\
&\phi(-\infty)=m_{\sigma}\quad\mbox{and}\quad\phi(+\infty)=m_{\sigma}+1,
\end{aligned}
\right.
\end{equation}

We start with the following statement about existence of traveling waves for $\sigma\in (\sigma^-,\sigma^+)$ and uniqueness of the velocity.
This is a reformulation of Theorem \ref{t2} i) above.

\begin{pro}\label{pro::r10}{\bf (Case $\sigma\in (\sigma^-,\sigma^+)$: uniqueness of the velocity; Theorem 1.2 and Theorem 1.6 a) in \cite{AFM})}\\
Assume that $F$ satisfies $(\tilde{A}_{C^1})$ and $(\tilde{B}_{C^1})$. For each value $\sigma\in (\sigma^-,\sigma^+)$, then there exists a velocity $c\in \R$ and a traveling wave $\phi$ solution of (\ref{e15bis}).
Moreover the velocity $c=c(\sigma)$ is unique.
\end{pro}

The remaining of this section is devoted to show some properties on the velocity function $c(\sigma)$.\\

In general we do not have sufficient regularity/assumptions in order to be able to compare the traveling waves (because we do not have strong comparison principle here).
Still we now show that we can compare the velocities.

\begin{pro}\label{p3}{\bf(Comparison of the velocities)}\\
Assume that $F$ satisfies $(\tilde{A}_{C^{1}})$ and $(\tilde{B}_{C^{1}})$. Let $\sigma_1,\sigma_2 \in [\sigma^-,\sigma^+]$ with
$$\sigma_1< \sigma_2$$ 
Assume that $(c_1,\phi_1)$ is a subsolution of (\ref{e15bis}) for $\sigma=\sigma_1$ and that $(c_2,\phi_2)$ is a supersolution of (\ref{e15bis}) for $\sigma=\sigma_2$ with
$$\phi_1(+\infty)> \phi_2(-\infty)$$
Then we have
$$c_{1}\leq c_{2}.$$
\end{pro}

\noindent {\bf Proof of Proposition \ref{p3}}\\
Recall that  $(c,\phi)$ is a solution of 
$$c\phi'(z)=F(\phi(z+r_{0}),\phi(z+r_{1}),...,\phi(z+r_{N}))+\sigma$$
if and only if
$$u(t,z):=\phi(z+ ct)$$
solves the evolution equation
$$u_t= F(u(z+r_{0}),u(z+r_{1}),...,u(z+r_{N}))+\sigma$$
The key point is that for this evolution equation, we always have a comparison principle under assumption $(\tilde{A}_{C^{1}})$ (see \cite[Propositions 2.5 and 2.6]{FIM}).\\
Now because of assumption $(\tilde{B}_{C^{1}})$, we have for the subsolution $\phi_1$ and the supersolution $\phi_2$
$$\left\{\begin{array}{l}
\phi_1(-\infty)\le m_{\sigma_1} < m_{\sigma_2}\le \phi_2(-\infty)\\
\phi_1(+\infty)\le 1+m_{\sigma_1} < 1+m_{\sigma_2}\le \phi_2(-\infty)\\
\end{array}\right.$$
and then up to translate $\phi_2$ for some $a\in \R$
$$\phi_2^a(z):=\phi_2(z+a)$$
we can insure
$$\phi_1 \le \phi_2^a$$
Setting
$$\left\{
\begin{aligned}
&u_{1}(t,z)=\phi_{1}(z+c_{1}t)\\
&u_{2}(t,z)=\phi_{2}^{a}(z+c_{2}t),
\end{aligned}
\right.
$$
we see that we can apply the comparison principle to the subsolution $u_1$ and the supersolution $u_2$ which are comparable at the initial time $t=0$.
This implies that
$$u_1(t,z)\le u_2(t,z)\quad \mbox{for all}\quad (t,z)\in [0,+\infty)\times \R$$
Then the fact that
$$\phi_1(+\infty)> \phi_2^a(-\infty)$$
implies the ordering of the velocities
$$c_1\le c_2$$
This ends the proof of the proposition.\\

\begin{cor}\label{c1}{\bf(Monotonicity of the velocities and the limits $c^\pm$)}\\
Assume $(\tilde{A}_{C^{1}}),$ $({\tilde B}_{C^{1}}).$
For $\sigma\in(\sigma^{-},\sigma^{+}),$ let $(c(\sigma),\phi_{\sigma})$ be a solution of (\ref{e15bis}) given in Proposition \ref{pro::r10}. Then the velocity function $\sigma \mapsto c(\sigma)$  is non-decreasing on $(\sigma^{-},\sigma^{+}).$ Moreover, the limits 
$$\lim_{\sigma^-<\sigma\to\sigma^{-}}c(\sigma)=c^{-}\quad\mbox{and}\quad\lim_{\sigma^+>\sigma\to\sigma^{+}}c(\sigma)=c^{+}$$ 
exist and satisfy $-\infty< c^{-}\leq c^{+}<+\infty.$
\end{cor}

\noindent {\bf Proof of Corollary \ref{c1}}\\
The monotonicity of the velocities follows immediately from Proposition \ref{p3}. This allows to define the limits $c^\pm$ which satisfy
$$-\infty \le c^{-}\leq c^{+} \le +\infty.$$
Let us show that $c^+$ is finite. From Proposition \ref{p0}, we know that for $\sigma=\sigma^+$, there exists a solution $(c_L,\phi)$ for $c_L>>1$ large enough. Then the comparison of the velocities (Proposition \ref{p3}) again gives for each $\sigma\in (\sigma^-,\sigma^+)$
$$c(\sigma)< c_L$$
This implies
$$c^+\le c_L<+\infty$$
and shows the finitness of $c^+$. The proof of the finitness of $c^-$ is similar. This ends the proof of the corollary.\\

In order to go further and show the continuity of the velocity function $\sigma\mapsto c(\sigma)$, we need a further ingredient that was also used  in the proof of Proposition \ref{pro::r10} to show the uniqueness of the velocity. This is the following comparison principle at infinity.

\begin{lem}\label{lem::r12}{\bf (Comparison principle at infinity; Theorem 4.1 and  Corrollary 4.2  in \cite{AFM})}\\
Assume that $F$  satisfies $(\tilde{A}_{C^{1}}),$ $({\tilde B}_{C^{1}})$ and let $\sigma\in (\sigma^-,\sigma^+)$ and any $c\in \R$.\\
\noindent {\bf i) (Comparison on $(-\infty,r^*]$)}\\
Let us consider $u,v: (-\infty,r^*] \to [m_\sigma,1+m_\sigma]$ be respectively a subsolution and a supersolution of the equation
$$cw'(z)=F((w(z+r_i))_{i=0,\dots,N})+\sigma \quad \mbox{on}\quad (-\infty,0)$$
Then there exists some $\delta_\sigma>0$ such that if
$$u\le \delta_\sigma+ m_\sigma \quad \mbox{on}\quad (-\infty,r^*]$$
then
$$u\le v \quad \mbox{on}\quad [0,r^*]$$
implies the comparison
$$u\le v \quad \mbox{on}\quad (-\infty,r^*]$$
\noindent {\bf ii) (Comparison on $[-r^*,+\infty)$)}\\
Let us consider $u,v: [-r^*,+\infty)\to [m_\sigma,1+m_\sigma]$ be respectively a subsolution and a supersolution of the equation
$$cw(z)=F((w(z+r_i))_{i=0,\dots,N})+\sigma \quad \mbox{on}\quad (0,+\infty)$$
Then there exists some $\delta_\sigma>0$ such that if
$$v\ge 1-\delta_\sigma+ m_\sigma \quad \mbox{on}\quad [-r^*,+\infty)$$
then
$$u\le v \quad \mbox{on}\quad [-r^*,0]$$
implies the comparison
$$u\le v \quad \mbox{on}\quad [-r^*,+\infty)$$
\end{lem}

\begin{pro}\label{p2}{\bf(Continuity of the velocity function)}\\
Suppose that $F$ satisfies $(\tilde{A}_{C^{1}}),$ $({\tilde B}_{C^{1}})$ and let $\sigma\in(\sigma^{-},\sigma^{+}).$ Let $(c(\sigma),\phi_{\sigma})$ be a solution of (\ref{e15bis}) given in Proposition \ref{pro::r10}. Then the map $\displaystyle\sigma\mapsto c(\sigma)$ is continuous on $(\sigma^{-},\sigma^{+}).$
\end{pro}

{\bf Proof of Proposition \ref{p2}}\\
Let $\sigma_{0}\in(\sigma^{-},\sigma^{+})$ and a sequence $\sigma_n\in (\sigma^{-},\sigma^{+})$, such that
$$\sigma_n \to \sigma_0\quad \mbox{as}\quad n\to +\infty$$
and let
$$c_0:=c(\sigma_0),\quad c_n:=c(\sigma_n)$$
Let $(c_0,\phi_0)$  and $(c_n,\phi_n)$ be  solutions of (\ref{e15bis}) respectively for the parameters $\sigma=\sigma_0$ and $\sigma=\sigma_n$, given by Proposition \ref{pro::r10}. Up to extract a subsequence, assume that
$$c_n\to c_\infty$$
We already know that $c_\infty \in [c^-,c^+]$ and want to show that $c_\infty= c_0$.\\
\noindent{\bf Case $1$: $c_0< c_\infty$}\\
Then up to shift $\phi_n$, we can assume that
$$\frac{m_{\sigma_0} + b_{\sigma_0}}{2} \in [\phi_n)_*(0),(\phi_n)^*(0)]$$
and then, up to extract a subsequence, we can assume that $\phi_n$ converges (at least almost everywhere) to some function $\phi_\infty$ which satisfies by stability (Proposition \ref{pro::r2})
$$\left\{\begin{array}{l}
\begin{aligned}
&c_\infty\phi_\infty'(z)=F((\phi_\infty(z+r_{i}))_{i=0,\dots,N})+\sigma_0\quad\mbox{on}\quad\R\\
&\phi_\infty\mbox{ is non-decreasing over }\R\\
&\phi_\infty(-\infty)\ge m_{\sigma_0}\quad\mbox{and}\quad\phi_\infty(+\infty)\le m_{\sigma_0}+1,\\
& \frac{m_{\sigma_0} + b_{\sigma_0}}{2} \in [\phi_\infty)_*(0),(\phi_\infty)^*(0)]
\end{aligned}
\end{array}\right.$$
which implies
$$\phi_\infty(-\infty)=m_{\sigma_0},\quad \phi_\infty(+\infty)\in \left\{b_{\sigma_0},1+m_{\sigma_0}\right\}$$
Hence $(c_0,\phi_\infty)$ is a subsolution of the equation satisfied by $(c_0,\phi_0)$ for $\sigma=\sigma_0$.
Then using the comparison at infinity (Lemma \ref{lem::r12}), we see that up to shift $\phi_0$, we can insure that
$$\phi_\infty \le \phi_0$$
Now, as in the proof of Proposition \ref{p3}, the comparison principle for the evolution equation shows that
$$\phi_\infty(x+ c_\infty t) \le \phi_0(x+c_0 t) \quad \mbox{for all}\quad (t,x)\in [0,+\infty)\times \R$$
which implies
$$c_\infty \le c_0$$
which is a contradiction.\\
\noindent{\bf Case $2$: $c_0> c_\infty$}\\
The reasoning is similar to the previous case. 
Here, up to shift $\phi_n$, we can now assume that
$$\frac{b_{\sigma_0}+1+m_{\sigma_0}}{2} \in [\phi_n)_*(0),(\phi_n)^*(0)]$$
and then, up to extract a subsequence, we can assume that $\phi_n$ converges (at least almost everywhere) to some function $\phi_\infty$ which satisfies by stability (Proposition \ref{pro::r2})
$$\left\{\begin{array}{l}
\begin{aligned}
&c_\infty\phi_\infty'(z)=F((\phi_\infty(z+r_{i}))_{i=0,\dots,N})+\sigma_0\quad\mbox{on}\quad\R\\
&\phi_\infty\mbox{ is non-decreasing over }\R\\
&\phi_\infty(-\infty)\ge m_{\sigma_0}\quad\mbox{and}\quad\phi_\infty(+\infty)\le m_{\sigma_0}+1,\\
& \frac{b_{\sigma_0}+1+m_{\sigma_0}}{2} \in [\phi_\infty)_*(0),(\phi_\infty)^*(0)]
\end{aligned}
\end{array}\right.$$
which implies
$$\phi_\infty(-\infty)\in \left\{m_{\sigma_0},b_{\sigma_0}\right\},\quad \phi_\infty(+\infty)=1+m_{\sigma_0}$$
Hence $(c_0,\phi_\infty)$ is a supersolution of the equation satisfied by $(c_0,\phi_0)$ for $\sigma=\sigma_0$.
Again using the comparison at infinity (Lemma \ref{lem::r12}), we see that up to shift $\phi_0$, we can insure that
$$\phi_\infty \ge \phi_0$$
Now, as in the proof of Proposition \ref{p3}, the comparison principle for the evolution equation shows that
$$\phi_\infty(x+ c_\infty t) \ge \phi_0(x+c_0 t) \quad \mbox{for all}\quad (t,x)\in [0,+\infty)\times \R$$
which implies
$$c_\infty \ge c_0$$
which is a contradiction.\\
\noindent{\bf Conclusion}\\
We deduce that $c_\infty=c_0$ and this shows the continuity of the velocity function. This ends the proof of the proposition.\\

\begin{lem}\label{l14}{\bf(Strict monotonicity)}\\
Assume that $F$ satisfies $(\tilde{A}_{C^1})$ and $(\tilde{B}_{C^1})$. Then  there exists a constant $K>0$ such that the velocity $c=c(\sigma)$ satisfies
\begin{equation}\label{e26}
\frac{dc}{d\sigma}\geq K|c|
\quad\mbox{on }\ (\sigma^{-},\sigma^{+})
\end{equation} 
in the viscosity sense.
\end{lem}

\noindent {\bf Proof of Lemma \ref{l14}}\\
Because of the  global Lipschitz continuity and periodicity of $F$ (by assumption $(\tilde{A}_{C^1})$), we know that there exists a constant $K>0$ such that we have the following useful bound
$$|F((\phi(z+r_i))_{i=0,\dots,N}) + \sigma| \le \frac{1}{K}$$
for any monotone function $\phi$ satisfying $\varphi(-\infty)=m_\sigma$ and $\varphi(+\infty)=1+m_\sigma$ uniformly in $\sigma \in (\sigma_-,\sigma_+)$.

Now let $\sigma_{1},\,\,\sigma_{2}\in(\sigma^{-},\sigma^{+})$ with $\sigma_{1}<\sigma_{2}$. Let $(c_1,\phi_1)$  and $(c_2,\phi_2)$ be  solutions of (\ref{e15bis}) respectively for the parameters $\sigma=\sigma_1$ and $\sigma=\sigma_2$, given by Proposition \ref{pro::r10}, with $c_i=c(\sigma_i)$ for $i=1,2$.\\
\noindent {\bf Case $c_1>0$}\\
Because we have
$$c_1 \phi_1'(x)= F((\phi_1(x+r_{i}))_{i=0,...,N})+\sigma$$
we deduce that
$$0\le \phi_1' \le \delta^{-1}\quad \mbox{with}\quad \delta:=K c_1$$
This implies that
$$\bar c \phi'_{1}\leq\sigma_{2}+F((\phi_{1}(x+r_{i}))_{i=0,...,N})\quad \mbox{with}\quad \bar c:= c_{1}+\delta(\sigma_{2}-\sigma_{1})$$
This means that $(\bar c,\phi_1)$ is a subsolution of the equation satisfied by $(c_2,\phi_2)$ for $\sigma=\sigma_2$. Up to shift $\phi_1$, we can assume that
$$\phi_1 \le \phi_2$$
and then, again, the comparison principle for the evolution equation implies that
$$\phi_1(x+\bar c t) \le \phi_2(x+c_2 t)\quad \mbox{for all}\quad (t,x)\in [0,+\infty)\times \R$$
We deduce that
$$\bar c \le c_2$$
ie
\begin{equation}\label{e30}
\frac{c_{2}-c_{1}}{\sigma_{2}-\sigma_{1}}\geq \delta = K c_1
\end{equation}
Now letting 
$\sigma_{2} \to \sigma_1$ or $\sigma_1 \to \sigma_2$, and using the continuity of $\sigma\mapsto c(\sigma)$, we deduce in the viscosity sense that
$$\frac{dc}{d \sigma} \ge K c\quad \mbox{on}\quad \left\{c>0\right\}\cap (\sigma^-,\sigma^+)$$
\noindent {\bf Case $c_1<0$}\\
We proceed similarly for $c_1<0$ and get the inequality on $\left\{c<0\right\}\cap (\sigma^-,\sigma^+)$.\\
\noindent {\bf Conclusion}\\
Finally, the monotonicity of $c$ implies that the desired inequality (\ref{e26}) holds true on the whole interval $(\sigma^-,\sigma^+)$. This ends the proof of the proposition.

\section{Vertical branches: filling the gaps and proof of Theorem \ref{t2}}\label{S5}

In this section, we focus on the study of solutions $(c,\phi)$ of (\ref{e17}), namely
\begin{equation}\label{e17ter}
\left\{
\begin{aligned}
&c\phi'(z)=F(\phi(z+r_{0}),\phi(z+r_{1}),...,\phi(z+r_{N}))+\sigma^{+}\quad\mbox{on}\quad\R\\
&\phi\mbox{ is non-decreasing over }\R\\
&\phi(-\infty)=0=m_{\sigma^{+}}\quad\mbox{and}\quad\phi(+\infty)=1.
\end{aligned}
\right.
\end{equation}
and the situation will be similar for $\sigma=\sigma^-$. Finally, at the end of this section we will be in position to give a proof of Theorem \ref{t2}.\\

\begin{lem}\label{l12}{\bf(Non-existence of solution for $c<c^{+}$ and $c>c^{-}$)}\\
Assume that $F$ satisfies $(\tilde{A}_{C^{1}})$ and $({\tilde B}_{C^{1}}).$ 
Let $(c,\phi)$ be a solution of (\ref{e17ter}) (resp. (\ref{e17'})), then $c\geq c^{+}$ (resp. $c\leq c^{-}$),
where $c^\pm$  are defined in Corollary \ref{c1}.
\end{lem}

\noindent {\bf Proof of Lemma \ref{l12}}\\
We do the proof for equation (\ref{e17ter})  (the proof is similar for (\ref{e17'})).
Let $(c,\phi)$ be a solution of \eqref{e17ter} for parameter $\sigma^+$.
From Proposition \ref{pro::r10}, it is known that for any $\sigma\in (\sigma^-,\sigma^+)$, there exists a solution $(c(\sigma),\phi_\sigma)$ of equation (\ref{e15bis}) with parameter $\sigma$.
Moreover the comparison of the velocities (Proposition \ref{p3}) implies that
$$c(\sigma) \le c$$
This implies that
$$c^+:=\lim_{\sigma^+>\sigma \to \sigma^+} c(\sigma) \le c$$
which implies the result. This ends the proof of the Lemma.\\

\begin{lem}\label{l2}{\bf(Existence of  solutions  for $c=c^{\pm}$)}\\
Assume $(\tilde{A}_{C^{1}}),$ $({\tilde B}_{C^{1}})$ and let $\sigma=\sigma^{+}$ (resp. $\sigma=\sigma^{-}$). There exists a profile $\phi^{+}$ (resp. $\phi^{-}$) such that $(c^{+},\phi^{+})$ (resp. $(c^{-},\phi^{-})$) solves (\ref{e17ter}) (resp. \eqref{e17'}), where $c^\pm$ are defined in Corollary \ref{c1}.
\end{lem}

\noindent {\bf Proof of Lemma \ref{l2}}\\
Assume that $\sigma=\sigma^{+}$ and let us prove the existence of a solution of (\ref{e17ter}) for $c^{+}$ (proving the existence of solution of (\ref{e17'}) for $c^{-}$ in the case $\sigma=\sigma^{-}$ is treated similarly). The goal is to get a solution as a limit of the profiles as $\sigma\to\sigma^{+},$ recalling that $\displaystyle c^{+}=\lim_{\sigma_+>\sigma\to\sigma^{+}}c(\sigma).$

Consider $\sigma\in(\sigma^{-},\sigma^{+})$ and let $(c_{\sigma},\phi_{\sigma})$ be a solution of (\ref{e15bis}) given by Proposition \ref{pro::r10}, namely 
\begin{equation}\label{e34}
\left\{
\begin{aligned}
&c_{\sigma}\phi'_{\sigma}(z)=F(\phi_{\sigma}(z+r_{0}),\phi_{\sigma}(z+r_{1}),...,\phi_{\sigma}(z+r_{N}))+\sigma\quad\mbox{on}\quad\R.\\
&\phi_{\sigma}\ \mbox{is non-decreasing over }\R\\
&\phi_{\sigma}(-\infty)=m_{\sigma}\quad\mbox{and}\quad\phi_{\sigma}(+\infty)=m_{\sigma}+1.
\end{aligned}
\right.
\end{equation} 
Up to shift $\phi_\sigma$, we can assume that
\begin{equation}\label{5e2}
\frac12 \in [(\phi_{\sigma})_{*}(0),(\phi_{\sigma})^{*}(0)]
\end{equation}
Passing to the limit $\sigma \to \sigma^+$, we get the convergence (at least almost everywhere) of $\phi_\sigma$ to some profile $\phi^+$.
Then by stability (Proposition \ref{pro::r2}), we get that $\phi^+$ solves
\begin{equation}\label{e35}
\left\{
\begin{aligned}
&c^{+}(\phi^+)'(z)=F((\phi^+(z+r_{i}))_{i=0,\dots,N})+\sigma^{+}\quad\mbox{on}\quad\R.\\
&\phi^+\ \mbox{is non-decreasing over }\R\\
&0=m_{\sigma^{+}}\leq\phi^+\leq m_{\sigma^{+}}+1=1\\
& \frac12\in [(\phi^+)_{*}(0),(\phi^+)^{*}(0)]
\end{aligned}
\right.
\end{equation}
which implies that
$$\phi^+(-\infty)=0,\quad \phi^+(+\infty)=1$$
This shows that $\phi^+$ is a desired profile with velocity $c^+$. This ends the proof.\\

In order to prove the existence of full vertical branches for $c\ge c^+$ (or $c\le c^-$) as in Proposition \ref{p4} below, it is natural to think to use Corollary \ref{sol->sup} which build a branch of solutions starting from a positive supersolution. This works well for instance if we know that for the velocity $c=c^+$, there  exists a positive profile $\phi^+>0$. The difficulty is that we have no strong maximum principle in our general setting, and then we can not assume that such profile is positive. Still we will construct a vertical branch of solutions for $c\ge c^+$, but we need the following result.

\begin{lem}\label{l3}{\bf(Existence of a hull function (\cite[Theorem 1.5 and Theorem 1.6, a1,a2]{FIM}))}\\ 
Assume that  $F$ satisfies  $(\tilde{A}_{C^{1}})$ and let $p>0$ and $\sigma\in\R.$ There exists a unique real $\lambda(\sigma,p)=\lambda_{p}(\sigma)$ such that there exists a locally bounded function $h_{p}:\R\to\R$ satisfying (in the viscosity sense):
$$\left\{
\begin{aligned}
&\lambda_{p}h'_{p}(z)=F((h_{p}(z+pr_{i}))_{i=0,...,N})+\sigma\quad\mbox{on}\quad\R\\&h_{p}(z+1)=h_{p}(z)+1\\
&h'_{p}(z)\geq 0\\
&|h_{p}(z+z')-h_{p}(z)-z'|\leq 1\quad\mbox{for any } z,\ z'\in\R.
\end{aligned}
\right.$$
Moreover, there exists a constant $K>0,$ independent on $p$ and $\sigma,$ such that 
$$|\lambda_{p}-\sigma|\leq K(1+p)$$
and the function 
$$
\begin{aligned}
\lambda_{p}:\ &\R\to\R\\
&\sigma\mapsto\lambda_{p}(\sigma)
\end{aligned}
$$
is continuous nondecreasing with $\lambda_{p}(\pm\infty)=\pm\infty.$
\end{lem}

Then we have
\begin{pro}\label{p4}{\bf(Existence of  vertical branches of velocities)}\\
Assume that  $F$ satisfies  $(\tilde{A}_{C^{1}})$ and $({\tilde B}_{C^{1}})$. Then for every velocity $c\ge c^{+}$ (resp. $c\le c^{-}$), there exists a solution $\phi$ of (\ref{e17ter}) (resp. (\ref{e17'})), where $c^\pm$ are defined in Corollary \ref{c1}.
\end{pro}

\noindent {\bf Proof of Proposition \ref{p4}}\\
\noindent {\bf Step 1: preliminary}\\
From the properties of $\lambda_p(\sigma)$, we deduce that for any velocity $c\in \R$, there exists a suitable $\sigma=\sigma(c,p)$ such that
$$\lambda_p(\sigma)=cp$$
Then defining the function $\phi_p$ as
$$\phi_{p}(x)=h_{p}(px),$$
where $h_{p}$, we see that it satisfies
\begin{equation}\label{e41}
\left\{
\begin{aligned}
&c\phi'_{p}(z)=F((\phi_{p}(z+r_{i}))_{i=0,...,N})+\sigma(c,p)\quad\mbox{on}\quad\R\\
&\phi'_{p}\ \mbox{non-decreasing}\\
&\phi_{p}\left(z+\frac{1}{p}\right)=\phi_{p}(z)+1.
\end{aligned}
\right.
\end{equation}
We will build a branch of solutions for $c\ge c^+$ (the proof is similar for $c\le c^-$).\\
\noindent{\bf Step 2: passing to the limit $p\to 0^+$}\\
Let us now consider the case of a fixed velocity
$$c\ge  c^+$$
Recall that in both cases $c=0$ and $c\not=0$, we know that $\phi_p$ is an almost everywhere solution of (\ref{e41}). Hence we can integrate it for instance on the interval $z\in [0,1]$, and from the periodicity and Lipschitz properties of $F$ deduce that
$$|\sigma(c,p)| \le C_0$$
where $C_0$ is independent on $p$. Up to translate the profile $\phi_p$, we can also assume that
$$\theta \in [(\phi_p)_*(0),(\phi_p)^*(0)]$$
Then we can pass to the limit $p\to 0$, and get the convergence $\phi_p \to \phi$, at least almost everywhere (and indeed everywhere if $c\not=0$) and the convergence
$$\sigma(c,p) \to \sigma$$
which satisfy
$$\left\{\begin{array}{l}
\begin{aligned}
&c\phi'(z)=F((\phi(z+r_{i}))_{i=0,...,N})+\sigma\\
&\phi\ \mbox{non-decreasing}\\
&\theta \in [(\phi)_*(0),(\phi)^*(0)]\\
&\phi(+\infty)-\phi(-\infty)\leq 1,
\end{aligned}
\end{array}\right.$$
Passing to the limit $z\to \pm \infty$, we get that
$$f(\phi(\pm \infty))+\sigma=0$$
which forces $\sigma\in [\sigma^-,\sigma^+]$ and
$$\phi(-\infty)=m_\sigma,\quad \phi(+\infty)=1+m_\sigma$$
\noindent{\bf Step 3: conclusion}\\
Assume by contradiction that
$$\sigma< \sigma^+$$
Then the comparison of velocities (Proposition \ref{p3}) implies that
$$c< c^+$$
which leads to a contradiction. We deduce that
$$\sigma=\sigma^+,\quad \phi(-\infty)=0,\quad \phi(+\infty)=1$$
and then $(c,\phi)$ is a solution of (\ref{e17ter}). This provides the existence of a profile for any velocity $c\ge c^+$.
This ends the proof of the proposition.\\

\begin{rem}\label{rem::r15}
Notice that the proof of Proposition \ref{p4} provides for the velocity $c=c^+$ a direct proof of the existence of profile $\phi^+$, different from the proof of Lemma \ref{l2}.
\end{rem}

As a corollary of the previous results, we now get the proof of Theorem \ref{t2}.\\

\noindent {\bf Proof of Theorem \ref{t2}}\\
\noindent {\bf Step 1: Proof of Theorem \ref{t2} 1) i)}\\
This part of the result follows from Proposition \ref{pro::r10} which provides existence of a solution $(c(\sigma),\phi_\sigma)$ with a unique velocity $c(\sigma)$ for  the range $\sigma\in  (\sigma^-,\sigma^+)$.\\
\noindent {\bf Step 2: Proof of Theorem \ref{t2} 1) ii)}\\
The bound from below 
$$\displaystyle \frac{d c}{d \sigma} \ge K c\quad \mbox{for}\quad \sigma\in (\sigma^-,\sigma^+)$$ 
follows from Lemma \ref{l14} on the strict monotonicity of the velocity.
The velocities $c^\pm$ are defined in Corollary \ref{c1} by
$$\lim_{\sigma^-<\sigma \to \sigma^-} c(\sigma)=:c^- \le c^+:=\lim_{\sigma^+>\sigma \to \sigma^+} c(\sigma)$$
We deduce that if $c^-\not=0$ of $c^-\not=0$, then we have
$$c^-<c^+$$
and then
$$c^-=c^+\quad \mbox{if only if}\quad c^-=0=c^+$$

\noindent {\bf Step 3: Proof of Theorem \ref{t2} 2) i) and ii)}\\
The proof of 2)ii) for $\sigma=\sigma^-$ is similar to the one of 2)i) for $\sigma=\sigma^+$. Hence we only show 2)i).
The non-existence of solutions for $c<c^+$ and $\sigma=\sigma^+$ follows from Lemma \ref{l12}. The existence of a branch of solutions for $c\ge c^+$ follows from Proposition \ref{p4}.\\
This ends the proof of the theorem.\\

\section{Appendix: example of discontinuous viscosity solutions}\label{Ap}

We give in this section an example of a discontinuous viscosity solution.

\begin{pro}\label{pdvs}{\bf(Discontinuous viscosity solution for the classical Frenkel-Kontorova model)}\\ 
Consider $\beta>0,$ $\sigma\in\R$ and let $(c,\phi)$ be a solution of 
\begin{equation}\label{2e1}
\left\{
\begin{aligned}
&c\phi'(z)=\phi(z+1)-2\phi(z)+\phi(z-1)+f(\phi(z))+\sigma\quad\mbox{on }\ \R,\quad f(x):=-\beta\cos (2\pi x) \\
&\phi\ \mbox{ is non-decreasing}\\
&\phi(+\infty)-\phi(-\infty)=1.
\end{aligned}
\right.
\end{equation}
\noindent {\bf i) (Sign of the critical velocities)}\\
Then $\sigma^{\pm}=\pm\beta$ and $c_-<0<c_+$.\\
\noindent {\bf ii) (Discontinuous solution for large $\beta$)}\\
Moreover, if $\beta>1$ and $|\sigma|<\beta-1,$ then $\phi\notin C^{0}$ and $c=0.$  
\end{pro}

\noindent

\noindent {\bf Proof of Proposition \ref{pdvs}}\\
\noindent {\bf Step 1: proof of i)}\\
Clearly, we have $\sigma^{\pm}=\pm\beta$. Let $\sigma=\sigma^{+}$ and let us show that $c^{+}>0.$ 
In this case, we can moreover assume that a solution $\phi$ of (\ref{2e1}) satisfies
$$\phi(-\infty)=0,\quad \phi(+\infty)=1$$ 
Integrating over the real line the equation 
$$c^{+}\phi'(z)=\phi(z+1)+\phi(z-1)-2\phi(z)+g(\phi(z)),\quad g:= f+ \sigma^+\ge 0$$ 
we get that 
$$c^{+}=\int_{\R}g(\phi(z))dz\geq0.$$ 
Since $g>0$ on $(0,1),$ if $c^{+}=0,$ then $$\phi(z)=0\mbox{ or }1\ \mbox{ almost everywhere}.$$ 
Then the equation itself implies that (because $g(\phi)=0$ a.e.)
\begin{equation}\label{eq::r1}
\Delta_{1}\phi(z):=\phi(z+1)+\phi(z-1)-2\phi(z)=0\ \mbox{ almost everywhere.}
\end{equation}
Because $\phi$ is monotone non-decreasing, up to translation, we have
$$\phi(z)=\left\{\begin{array}{ll}
0&\quad \quad \mbox{if}\quad z<0\\
1&\quad \quad \mbox{if}\quad z>0\\
\end{array}\right.$$
This leads to a contradiction with (\ref{eq::r1}), and shows that $c_+>0$. Similarly, we get $c_-<0$.\\
\noindent {\bf Step 2: proof of ii)}\\
Let $|\sigma|<\beta-1$ and let us show that $\phi\notin C^{0}(\R).$ 
For the convenience of the reader we give the proof of this result (which is basically contained in Theorem $1.2$ in Carpio et al. \cite{CCHM}).\\
Assume to the contrary that $\phi\in C^{0}(\R)$.
Notice that because $\phi$ is non-decreasing and $\phi(+\infty)-\phi(-\infty)=1,$ we deduce that $$\phi(z+1)-2\phi(z)+\phi(z-1)\in[-1,1].$$
Define now 
$$\psi(z)=\phi(z+1)-2\phi(z)+\phi(z-1)+f(\phi(z))+\sigma \quad \mbox{with}\quad f(\phi):=-\beta \cos (2\pi \phi)$$ 
Assume by contradiction that $\phi\in C^{0}$. Then we deduce that 
$$\left\{
\begin{aligned}
&\sup_{\R}\psi\geq\beta+\sigma-1>0\\
&\inf_{\R}\psi\leq-\beta+\sigma+1<0,
\end{aligned}
\right.$$
where the strict inequalities follow from $|\sigma|<\beta-1.$
But $c\phi'=\psi$ which implies that $c\phi'$ changes sign. Contradiction. Therefore  $\phi\notin C^{0}(\R)$, which also implies that $c=0.$\\
This ends the proof of the proposition.\\

\vspace{10mm}

\noindent{\bf Acknowledgments}

The first author would like to thank the Lebanese National Council for Scientific Research (CNRS-L) and the Campus France (EGIDE earlier) for supporting him. He also want to thank professor R. Talhouk and the Lebanese university. The last author was also partially supported by the contract ERC ReaDi 321186. Finally, this work was partially supported by ANR HJNet (ANR-12-BS01-0008-01) and by ANR-12-BLAN-WKBHJ: Weak KAM beyond Hamilton-Jacobi.\\

There are no other funding. This work has been done in collaboration. There is no conflict of interest.


\end{document}